\documentclass{amsart}
\usepackage{amscd}
\usepackage{amssymb}

\def\Kweb#1{
http:\linebreak[3]//www.\linebreak[3]math.\linebreak[3]uiuc.
\linebreak[3]edu/\linebreak[3]{K-theory/#1/}}

\newcommand{\cA}{\mathcal{A}}
\newcommand{\cB}{\mathcal{B}}
\def\cC{\mathcal C}

\def\EE{\mathcal E}
\def\cF{\mathcal F}

\def\cO{\mathcal O}
\def\cP{\mathcal P}
\def\cK{\mathcal K}
\def\cKH{\mathcal{KH}}
\def\bfP{\mathbf P}
\newcommand{\A}{\mathbb{A}}
\newcommand{\C}{\mathbb{C}}
\newcommand{\Q}{\mathbb{Q}}

\newcommand{\Z}{\mathbb{Z}}
\newcommand{\bP}{\mathbb{P}}

\def\bu{\bullet}

\def\Hom{\operatorname{Hom}}

\def\Spec{\operatorname{Spec}}

\newcommand{\SchF}{\mathrm{Sch}/F}
\newcommand{\eps}{\varepsilon}
\newcommand{\SmF}{\mathrm{Sm}/F}

\newcommand{\D}{\operatorname{D_{parf}}}
\newcommand{\oD}{\operatorname{D}}
\newcommand{\Df}[1]{\operatorname{D}^{#1}_{\mathrm{parf}}}
\newcommand{\bCh}{\operatorname{\mathbf{Ch}}}
\newcommand{\Ac}{\operatorname{Ac}}
\newcommand{\cS}{\mathcal{S}}
\newcommand{\cT}{\mathcal{T}}
\newcommand{\dL}{\mathrm{L}}
\newcommand{\dR}{\mathrm{R}}
\newcommand{\dLq}{\dL q}
\newcommand{\dRq}{\dR q}
\newcommand{\dLj}{\dL j}
\newcommand{\dRj}{\dR j}
\newcommand{\dLp}{\dL p}
\newcommand{\dRp}{\dR p}

\def\bbH{\mathbb H}
\newcommand{\bHN}{\mathbf{HN}}
\newcommand{\bHP}{\mathbf{HP}}
\newcommand{\bHC}{\mathbf{HC}}
\newcommand{\Mix}{\mathcal Mix}
\newcommand{\DMix}{\mathcal{D}\Mix}
\def\bC{\mathbf C}

\input xypic
\xyoption{all} \CompileMatrices

\numberwithin{equation}{section}

\theoremstyle{plain}
\newtheorem{thm}[equation]{Theorem}
\newtheorem{cor}[equation]{Corollary}
\newtheorem{lem}[equation]{Lemma}
\newtheorem{prop}[equation]{Proposition}
\newtheorem{conj}[equation]{$K$-dimension Conjecture}

\theoremstyle{definition}
\newtheorem{defn}[equation]{Definition}
\theoremstyle{remark}
\newtheorem{rem}[equation]{Remark}
\newtheorem{ex}[equation]{Example}

\newtheorem{terminology}[equation]{Terminology}
\newtheorem{Exact-cat}[equation]{Exact categories}
\newtheorem{Exactdg}[equation]{Exact dg categories}
\newtheorem{locpairs}[equation]{Localization pairs}
\newtheorem{sqlocpairs}[equation]{Sub and quotient localization pairs}

\begin{document}
\bibliographystyle{plain}

\title{Cyclic homology, $cdh$-cohomology and negative $K$-theory}

\author{G. Corti\~nas}
\thanks{Corti\~nas' research was partially supported by the Ram\'on
y Cajal fellowship, by ANPCyT grant PICT 03-12330 and by MEC grant
MTM00958.}
\address{Dep.\ \'Agebra y Geometr\'\i a y Topolog\'\i a,
Universidad de Valladolid, 47005 Valladolid, Spain}
\email{gcorti@agt.uva.es}

\author{C. Haesemeyer}
\thanks{Haesemeyer's research was partially supported by the Bell Companies
Fellowship and RTN Network HPRN-CT-2002-00287.}
\address{Dept.\ of Mathematics, University of Illinois, Urbana, IL
61801} \email{chh@math.uiuc.edu}

\author{M. Schlichting}
\thanks{Schlichting's research was partially supported by
RTN Network HPRN-CT-2002-00287.}
\address{Dept.\ of Mathematics, Louisiana State University, Baton
Rouge, LA 70803} \email{mschlich@math.lsu.edu}

\author{C. Weibel}
\thanks{Weibel's research was partially supported by NSA grant MSPF-04G-184.}
\address{Dept.\ of Mathematics, Rutgers University, New Brunswick,
NJ 08901} \email{weibel@math.rutgers.edu}

\date{\today}

\begin{abstract}
We prove a blow-up formula for cyclic homology which we use to
show that infinitesimal $K$-theory satisfies $cdh$-descent.
Combining that result with some computations of the
$cdh$-cohomology of the sheaf of regular functions, we verify a
conjecture of Weibel predicting the vanishing of algebraic
$K$-theory of a scheme in degrees less than minus the dimension of
the scheme, for schemes essentially of finite type over a field of
characteristic zero.
\end{abstract}

\maketitle

\section*{Introduction}

The negative algebraic $K$-theory of a singular variety is related
to its geometry. This observation goes back to the classic study
by Bass and Murthy \cite{BM}, which implicitly calculated the
negative $K$-theory of a curve $X$. By definition, the group
$K_{-n}(X)$ describes a subgroup of the Grothendieck group $K_0(Y)$
of vector bundles on $Y=X\times(\A^1-\{0\})^n$.

The following conjecture was made in 1980, based upon the Bass-Murthy
calculations, and appeared in \cite[2.9]{WeibelKA}. Recall that if $F$
is any contravariant functor on schemes, a scheme $X$ is called
$F$-regular if $F(X)\to F(X\times\A^r)$ is an isomorphism for all
$r\ge0$.

\begin{conj}\label{conj:reg}
Let $X$ be a Noetherian scheme of dimension $d$. Then $K_m(X) = 0$
for $m < -d$ and $X$ is $K_{-d}$-regular.
\end{conj}

In this paper we give a proof of this conjecture for $X$
essentially of finite type over a field $F$ of characteristic $0$;
see Theorem \ref{thm:Kdim}. We remark that this conjecture is
still open in characteristic $p>0$, except for curves and
surfaces; see \cite{WeibelNorm}.

Much of this paper involves cohomology with respect to Voevodsky's
$cdh$-topology. The following statement summarizes some of our
results in this direction:

\begin{thm}\label{thm:cdh-intro}
Let $F$ be a field of characteristic $0$, $X$ a $d$-dimensional scheme,
essentially of finite type over $F$. Then:
\begin{enumerate}
\item $K_{-d}(X)\cong H^d_{cdh}(X,\Z)$ (see \ref{thm:Kdim});
\item $H^d_{Zar}(X,\cO_X) \to H^d_{cdh}(X,\cO_X)$ is surjective
(see \ref{thm:topchange});
\item If $X$ is smooth then $H^n_{Zar}(X,\cO_X) \cong H^n_{cdh}(X,\cO_X)$
for all $n$ (see \ref{prop:smoothtc}).
\end{enumerate}
\end{thm}

In addition to our use of the $cdh$-topology, our key technical
innovation is the use of Corti\~nas' infinitesimal $K$-theory
\cite{CortInf} to interpolate between $K$-theory and cyclic
homology. We prove (in Theorem \ref{thm:Kinfdesc}) that
infinitesimal $K$-theory satisfies descent for the $cdh$-topology.
Since we are in characteristic zero, every scheme is locally
smooth for the $cdh$-topology, and therefore locally $K_n$-regular
for every $n$. In addition, periodic cyclic homology is locally de
Rham cohomology in the $cdh$-topology. These features allow us to
deduce conjecture \ref{conj:reg} from Theorem \ref{thm:cdh-intro}.

This paper is organized as follows. The first two sections study
the behavior of cyclic homology and its variants under blow-ups.
We then recall some elementary facts about descent for the
$cdh$-topology in section~3, and provide some examples of functors
satisfying $cdh$-descent, like periodic cyclic homology
(\ref{cor:HPdesc}) and homotopy $K$-theory (\ref{ex:KHdesc}). We
introduce infinitesimal $K$-theory in section~4 and prove that it
satisfies $cdh$-descent. This already suffices to prove that $X$
is $K_{-d-1}$-regular and $K_n(X) = 0$ for $n< -d$, as
demonstrated in section~5. The remaining step, involving $K_{-d}$,
requires an analysis of the $cdh$-cohomology of the structure
sheaf $\cO_X$ and is carried out in section~6.
\\
\smallskip

{\bf Notation}

The category of spectra we use in this paper will not be critical.
In order to minimize technical issues, we will use the terminology
that a {\it spectrum} $E$ is a sequence $E_n$ of simplicial sets
together with bonding maps $b_n: E_n \to \Omega E_{n+1}$. We say
that $E$ is an {\it $\Omega$-spectrum} if all bonding maps are
weak equivalences. A map of spectra is a strict map. We will use
the model structure on the category of spectra defined in
\cite{BF}. Note that in this model structure, every fibrant
spectrum is an $\Omega$-spectrum.

If $A$ is a ring, $I\subset A$ a two-sided ideal and $\EE$ a
functor from rings to spectra, we write $\EE(A,I)$ for the
homotopy fiber of $\EE(A) \to \EE(A/I)$. If moreover $f:A \to B$
is a ring homomorphism mapping $I$ isomorphically to a two-sided
ideal (also called $I$) of $B$, then we write $\EE(A,B,I)$ for the
homotopy fiber of the natural map $\EE(A,I) \to \EE(B,I)$. We say
that $\EE$ {\it satisfies excision} provided that
$\EE(A,B,I)\simeq 0$ for all $A$, $I$ and $f: A\to B$ as above. Of
course, if $\EE$ is only defined on a smaller category of rings,
such as commutative $F$-algebras of finite type, then these
notions still make sense and we say that $\EE$ satisfies excision
for that category.

We shall write $\SchF$ for the category of schemes essentially of
finite type over a field $F$. We say a presheaf $\EE$ of spectra
on $\mathrm{Sch}/F$ {\it satisfies the Mayer-Vietoris-property}
(or MV-property, for short) for a cartesian square of schemes
$$
\begin{CD}
Y' @>>> X' \\
@VVV @VVV \\
Y @>>> X
\end{CD}
$$
if applying $\EE$ to this square results in a homotopy cartesian
square of spectra.  We say that $\EE$ satisfies the Mayer-Vietoris
property for a class of squares provided it satisfies the
MV-property for each square in the class. For example, the
MV-property for affine squares in which $Y\to X$ is a closed
immersion is the same as the excision property for commutative
algebras of finite type, combined with invariance under
infinitesimal extensions.

We say that $\EE$ {\it satisfies Nisnevich descent} for
$\mathrm{Sch}/F$ if $\EE$ satisfies the MV-property for all elementary
Nisnevich squares in $\mathrm{Sch}/F$; an {\em elementary Nisnevich
square} is a cartesian square of schemes as above for which $Y\to X$
is an open embedding, $X'\to X$ is \'etale and $(X'-Y')\to(X-Y)$ is an
isomorphism.  By \cite[4.4]{Nis}, this is equivalent to the assertion
that $\EE(X) \to \bbH_{nis}(X,\EE)$ is a weak equivalence for
each scheme $X$, where $\bbH_{nis}(-,\EE)$ is a fibrant
replacement for the presheaf $\EE$ in a suitable model structure.

We say that $\EE$ {\it satisfies $cdh$-descent} for
$\mathrm{Sch}/F$ if $\EE$ satisfies the MV-property for all
elementary Nisnevich squares (Nisnevich descent) and for all
abstract blow-up squares in $\mathrm{Sch}/F$. Here an {\em
abstract blow-up square} is a square as above such that $Y\to X$
is a closed embedding, $X'\to X$ is proper and the induced
morphism $(X' - Y')_{red} \to (X-Y)_{red}$ is an isomorphism. We
will see in Theorem \ref{thm:BG} that this is equivalent to the
assertion that $\EE(X) \to \bbH_{cdh}(X,\EE)$ is a weak
equivalence for each scheme $X$, where $\bbH_{cdh}(-,\EE)$ is a
fibrant replacement for the presheaf $\EE$ in a suitable model
structure.

It is well known that there is an Eilenberg-Mac~Lane functor from
chain complexes of abelian groups to spectra, and from presheaves
of chain complexes of abelian groups to presheaves of spectra.
This functor sends quasi-isomorphisms of complexes to weak
homotopy equivalences of spectra. In this spirit, we will use the
above descent terminology for presheaves of complexes. Because we
will eventually be interested in hypercohomology, we use
cohomological indexing for all complexes in this paper; in
particular, for a complex $A$, $A[p]^q = A^{p+q}$.

\section{Perfect complexes and regular blowups}

In this section, we compute the categories of perfect complexes
for blow-ups along regularly embedded centers. Our computation
slightly differs from that of Thomason (\cite{ThomRD}) in that we
use a different filtration which is more useful for our purposes.
We don't claim much originality.

In this section, ``scheme'' means ``quasi-separated and
quasi-compact scheme''. For such a scheme $X$, we write
  $\D(X)$ for the derived category of perfect complexes on $X$ \cite{TT}.
Let $i: Y \subset X$ be a regular embedding of schemes of pure
codimension $d$, and let $p: X' \to X$ be the blow-up of $X$ along
$Y$ and $j: Y' \subset X'$ the exceptional divisor. We write $q$
for the map $Y' \to Y$.

Recall that the exact sequence of $\cO_{X'}$-modules $\cO_{X'}(1)
\to \cO_{X'} \to j_*\cO_{Y'}$ gives rise to the fundamental exact
triangle in $\D(X')$:
\begin{equation}\label{M:eqn1}
\cO_{X'}(l+1) \to \cO_{X'}(l) \to \dRj_*\bigl(\cO_{Y'}(l)\bigr)
\to \cO_{X'}(l+1)[1],
\end{equation}
where $\dRj_*\bigl(\cO_{Y'}(l)\bigr) =\bigl(j_*\cO_{Y'}\bigr)(l)$
by the projection formula.

We say that a triangulated subcategory $\cS\subset \cT$ of a
triangulated category $\cT$ is generated by a specified set of
objects of $\cT$ if $\cS$ is the smallest {\it thick} (that is,
closed under direct factors) triangulated subcategory of $\cT$
containing that set.

\begin{lem}
\label{M:lemma2}
\begin{enumerate}
\item The triangulated category $\D(X')$ is generated by
$\dLp^*F$, $\dRj_* \dLq^* G \otimes \cO_{X'}(-l)$, for $F \in
\D(X)$, $G \in \D(Y)$ and $l = 1,..., d-1$. \item The triangulated
category $\D(Y')$ is generated by $\dLq^* G \otimes \cO_{Y'}(-l)$,
for $G \in \D(Y)$ and $l = 0,..., d-1$.
\end{enumerate}
\end{lem}

\begin{proof}
(Thomason \cite{ThomRD}) For $k=0,...,d$, let $\cA'_k$ denote the
full triangulated subcategory of $\D(X')$ of those complexes $E$
for which $\dRp_*(E \otimes \cO_{X'}(l))=0$ for $0 \leq l < k$. In
particular, $\D(X')=\cA'_0$. By \cite[Lemme 2.5(b)]{ThomRD},
$\cA'_d=0$. Using \cite[Lemme 2.4(a)]{ThomRD}, and descending
induction on $k$, we see that for $k \geq 1$, $\cA'_k$ is
generated by $\dRj_*\dLq^*G \otimes \cO_{X'}(-l)$, for some $G$ in
$\D(Y)$ and $l=k,...,d-1$. For $k=0$, we use the fact that the
unit map $1 \to \dRp_*\dLp^*$ is an isomorphism \cite[Lemme
2.3(a)]{ThomRD} to see that $\cA'_0=\D(X')$ is generated by the
image of $\dLp^*$ and the kernel of $\dRp_*$. But $\cA'_1$ is the
kernel of $\dRp_*$.

Similarly, for $k=0,...,d$, let $\cA_k$ be the full triangulated
subcategory of $\D(Y')$ of those complexes $E$ for which
$\dRq_*(E\otimes\cO_{Y'}(l))=0$ for $0 \leq l < k$. In particular,
$\D(Y')=\cA_0$. By \cite[Lemme 2.5(a)]{ThomRD}, $\cA'_d=0$. Using
\cite[p.247, from ``Soit $F^{\cdot}$ un objet dans
$\underline{\cA'}_k$'' to ``Alors $G^{\cdot}$ est un objet dans
$\underline{\cA'}_{k+1}$'']{Tproj}, and descending induction on
$k$, we have that $\cA_k$ is generated by $\dLq^*G
\otimes\cO_{Y'}(-l)$, $l=k,...,d-1$.
\end{proof}

\begin{rem}\label{M:rem1}
As a consequence of the proof of \ref{M:lemma2}, we
note the following. Let $k=0,...,d-1$ and $m$ any integer. The
full triangulated subcategory of $\D(Y')$ of those complexes $E$
with $\dRq_*(E \otimes \cO_{Y'}(l))=0$ for $m \leq l < k+m$ is
the same as the full triangulated subcategory
generated by $\dLq^*G \otimes \cO_{Y'}(n)$,
for $G \in \D(Y)$ and $k+m \leq n \leq d-1+m$.
In particular, the condition that a complex be in the latter category
is local in $Y$.
\end{rem}

\begin{lem}\label{M:lemma1}
The functors $\dLp^*\!:\! \D(X)\!\to\!\D(X')$,
$\dLq^*:\! \D(Y)\!\to\!\D(Y')$ and
$\dRj_* \dLq^*\!:\! \D(Y)\!\to\!\D(X')$
are fully faithful.
\end{lem}

\begin{proof}
The functors $\dLp^*$ and $\dLq^*$
are fully faithful, since the unit maps $1 \to \dRp_*\dLp^*$ and
$1 \to \dRq_*\dLq^*$ are isomorphisms \cite[Lemme 2.3]{ThomRD}.

By the fundamental exact triangle (\ref{M:eqn1}), the cone of the
counit $\dLj^* \dRj_* \cO_{Y'} \to \cO_{Y'}$ is in the
triangulated subcategory generated by $\cO_{Y'}(1)$, since the
counit map is a retraction of $\dLj^*\cO_{X'} \to
\dLj^*\dRj_*\cO_{Y'}$. It follows that the cone of the counit map
$\dLj^* \dRj_* \dLq^* E \to \dLq^* E$ is in the triangulated
subcategory generated by $\dLq^* E \otimes \cO_{Y'}(1)$, since the
latter condition is local in $Y$ (see Remark \ref{M:rem1}), and
$\D(Y)$ is generated by $\cO_Y$ for affine $Y$.  Since
$\dRq_*(\dLq^*G \otimes \cO(-1)) = G \otimes \dRq_* \cO(-1) = 0$,
we have $\Hom(A,B) = 0$ for $A$ (respectively $B$) in the
triangulated subcategory of $\D(Y')$ generated by $\dLq^*G \otimes
\cO(1)$ (respectively, generated by $\dLq^*G$), for $G \in \D(Y)$.
Applying this observation to the cone of
$\dLj^*\dRj_*\dLq^*E\to\dLq^*E$ justifies the second equality in
the display:
\begin{eqnarray*}
\Hom(E,F) = \Hom(\dLq^* E, \dLq^* F)
 =& \Hom(\dLj^* \dRj_* \dLq^* E, \dLq^* F) \\
 =& \Hom(\dRj_* \dLq^* E, \dRj_*\dLq^*F).
\end{eqnarray*}
The first equality holds because $\dLq^*$ is fully faithful, and
the final equality is an adjunction. The composition is an
equality, showing that $\dRj_*\dLq^*$ is fully faithful.
\end{proof}

For $l = 0,..., d-1$, let $\Df{l}(X') \subset\D(X')$ be the full
triangulated subcategory generated by $\dLp^*F$ and $\dRj_*\dLq^*G
\otimes \cO_{X'}(-k)$ for $F \in \D(X)$, $G \in \D(Y)$ and $k =
1,..., l$. For $l = 0,..., d-1$, let
$\oD^l_{\mathrm{parf}}(Y')\subset\D(Y')$ be the full triangulated
subcategory generated by $\dLq^*G \otimes \cO_{Y'}(-k)$ for $G \in
D(Y)$ and $k = 0,..., l$.

By Lemma \ref{M:lemma1}, $\dLp^*: \D(X) \to
\oD^0_{\mathrm{parf}}(X')$ and $\dLq^*: \D(Y) \to
\oD^0_{\mathrm{parf}}(Y')$ are equivalences. By Lemma
\ref{M:lemma2}, $\oD^{d-1}_{\mathrm{parf}}(X') = \D(X')$ and $
\oD_{\mathrm{parf}}^{d-1}(Y) = \D(Y')$.

\begin{prop}
\label{M:prop:filt}
The functor $\dLj^*$ is compatible with
the filtrations on $\D(X')$ and $\D(Y')$:
$$\xymatrix{\D(X) \ar[r]^{\dLp^*}_{\sim}
 \ar[d]^{\dL i^*} & \oD^0_{\mathrm{parf}}(X') \ar[d]^{\dLj^*}
 \ar@{}[r]|{\subset} & \oD^1_{\mathrm{parf}}(X')
 \ar@{}[rr]|<<<<<<<<<<{\subset\ \  \cdots\ \  \subset} \ar[d]^{\dLj^*} &&
\oD_{\mathrm{parf}}^{d-1}(X') = \D(X') \ar[d]^{\dLj^*} \\
\D(Y) \ar[r]^{\dLq^*}_{\sim}  & \oD^0_{\mathrm{parf}}(Y')
 \ar@{}[r]|{\subset} & \oD^1_{\mathrm{parf}}(Y')
 \ar@{}[rr]|<<<<<<<<<<{\subset \ \ \cdots\ \  \subset} &&
 \oD_{\mathrm{parf}}^{d-1}(Y') = \D(Y').
}$$
For $l=0,...,d-2$, $\dLj^*$ induces equivalences on successive
quotient triangulated categories:
$$\dLj^*: \Df{l+1}(X')/\Df{l}(X') \stackrel{\simeq}{\longrightarrow}
           \Df{l+1}(Y')/\Df{l}(Y').
$$
\end{prop}

\begin{proof}
The commutativity of the left hand square follows from $\dLq^* \dL
i^* = \dLj^* \dLp^*$. The compatibility of $\dLj^*$ with the
filtrations only needs to be checked on generators, that is, we
need to check that  $\dLj^* [\dRj_*\dLq^*G \otimes \cO_{X'}(-l)]$
is in $\Df{l}(Y')$, $l=1,...,d-1$. The last condition is local in
$Y$ (see Remark \ref{M:rem1}), a fortiori it is local in $X$. So
we can assume that $X$ and $Y$ are affine, and $G=\cO_Y$. In this
case, the claim follows from the fundamental exact triangle
(\ref{M:eqn1}).

For $l-k =1,...,d-1$, $E \in \D(X)$ and $G \in \D(Y)$, we have
$\Hom(\dLp^* E \otimes \cO(-k), \dRj_*\dLq^*G \otimes \cO(-l)) =
\Hom(\dLp^* E \otimes \cO(l-k), \dRj_*\dLq^*G) = \Hom(\dLj^*\dLp^*
E \otimes\cO(l-k), \dLq^*G) = \Hom(\dLq^*\dL i^* E \otimes
\cO(l-k), \dLq^*G) = 0$ since $\dRq_* \cO(k-l) =0$. Therefore, all
maps from objects of $\Df{l}(X')$ to an object of $\cO(-l-1)
\otimes\dRj_*\dLq^* \D(Y) \subset \Df{l+1}(X')$ are trivial. It
follows that the composition
$$\cO(-l-1) \otimes \dRj_*\dLq^* \D(Y) \subset
  \Df{l+1}(X') \to \Df{l+1}(X') / \Df{l}(X')
$$
is an equivalence (it is fully faithful, both categories have the
same set of generators, and the source category is idempotent
complete).
Similarly, the composition
$$\cO(-l-1) \otimes \dLq^* \D(Y) \subset \Df{l+1}(Y') \to
\oD^{l+1}_{\mathrm{parf}}(Y') / \oD^{l}_{\mathrm{parf}}(Y')$$ is
an equivalence.

The counit map $\dLj^*\dRj_*\dLq^* \to \dLq^*$ has cone in the
triangulated subcategory  generated by $\dLq^*G \otimes \cO(1)$
(see proof of \ref{M:lemma1}), $G \in \D(Y)$. It follows that the
natural map of functors $\dLj^*[\cO(-l-1) \otimes \dRj_*\dLq^*]
\to \cO_{Y'}(-l-1) \otimes \dLq^*$, induced by the counit of
adjunction, has cone in $\oD^l_{\mathrm{parf}}(Y')$. Thus, the
composition $\dLj^* \circ [\cO(-l-1) \otimes \dRj_*\dLq^*]: \D(Y)
\to \Df{l+1}(X')/\Df{l}(X') \to \Df{l+1}(Y')/\Df{l}(Y')$ agrees,
up to natural equivalence of functors, with $\cO_{Y'}(-l-1)
\otimes\dLq^*: \D(Y) \to \Df{l+1}(Y')/\Df{l}(Y')$. Since two of
the three functors are equivalences, so is the third: $\dLj^*:
\Df{l+1}(X')/\Df{l}(X') \stackrel{\sim}{\to}
\Df{l+1}(Y')/\Df{l}(Y')$.
\end{proof}

\begin{rem}
\label{M:rem2} Proposition \ref{M:prop:filt} yields $K$-theory
descent for blow-ups along regularly embedded centers.
This follows from Thomason's theorem in \cite{TT} (see \cite{GilletRR,GS}),
because every square in \ref{M:prop:filt}
induces a homotopy cartesian square of $K$-theory spectra.

Several people have remarked that this descent also follows from the
main theorem of \cite{ThomRD} by a simple manipulation.
\end{rem}

\section{Thomason's theorem for (negative) cyclic homology}

In this section we prove that negative cyclic, periodic cyclic and
cyclic homology satisfy the Mayer-Vietoris property for blow-ups
along regularly embedded centers. We will work over a ground field
$k$, so that all schemes are $k$-schemes, all linear categories
are $k$-linear, and tensor product $\otimes$ means tensor product
over $k$.

\subsection*{Mixed Complexes}

In order to fix our notation, we recall some standard definitions
(see \cite{LodayHC92} and \cite{WeibelHA94}). We remind the reader
that we are using cohomological notation, with the homology of $C$
being given by $H_n(C)=H^{-n}(C)$.

A {\it mixed complex} $\mathbf{C}=(C,b,B)$ is a cochain complex
$(C,b)$ ($C^n=0$ for $n>0$ is typical), together with a chain map
$B:C\to C[-1]$ satisfying $B^2=0$.
There is an evident notion of a map of mixed complexes, and we
write $\Mix$ for the category of mixed complexes.

The complexes for cyclic, periodic cyclic and negative cyclic
homology of $(C,b,B)$ are obtained using the product total
complex:
$$
\begin{array}{rclcccccccccc}
HC(C,b,B)&=&{\rm Tot} (\cdots\to  &
C[+1] &\stackrel{B}{\to} &  C & {\to}& 0 & {\to}& 0 & \to\cdots\ )\\
HP(C,b,B)&=&{\rm Tot} (\cdots\to & C[+1] &\stackrel{B}{\to} &  C &
\stackrel{B}{\to}& C[-1] &
            \stackrel{B}{\to} &C[-2] & \to\cdots\ )\\
HN(C,b,B)&=&{\rm Tot} (\cdots\to & 0 &{\to} &  C &
\stackrel{B}{\to} &C[-1] & \stackrel{B}{\to} &C[-2] & \to\cdots\ )
\end{array}
$$
where $C$ is placed in horizontal degree $0$. In particular, the
familiar exact sequence $0\to C\to HC(\mathbf{C})\to
HC(\mathbf{C})[+2]\to0$ is paired with a natural exact sequence of
complexes
$$0 \to HN(\mathbf{C}) \to HP(\mathbf{C}) \to HC(\mathbf{C})[+2] \to 0.$$
Short exact sequences and quasi-isomorphisms of mixed complexes
yield short exact sequences and quasi-isomorphisms of $HC$, $HP$
and $HN$ complexes, respectively. Of course, the cyclic, periodic
cyclic and negative cyclic homology groups of $\mathbf{C}$ are the
homology groups of $HC$, $HP$ and $HN$, respectively.

We say that a map $(C,b,B)\to(C',b',B')$ is a quasi-isomorphism in
$\Mix$ if the underlying complexes are quasi-isomorphic via
$(C,b)\to(C',b')$; following \cite{Keller99}, we write $\DMix$ for
the localization of $\Mix$ with respect to quasi-isomorphisms; it
is a triangulated category with shift
$\mathbf{C}\mapsto\mathbf{C}[1]$. The reader should beware that
$\DMix$ is {\it not} the derived category of the underlying
abelian category of $\Mix$.

It is sometimes useful to use the equivalence between the category
$\Mix$ of mixed complexes and the category of left dg
$\Lambda$-modules, where $\Lambda$ is the dg-algebra
$$\cdots 0 \to  k \eps \stackrel{0}{\to} k \to 0 \to \cdots$$
with $k$ placed in degree zero \cite[2.2]{Keller98}.
A left dg $\Lambda$ module $(C,d)$ corresponds to the mixed
complex $(C,b,B)$ with $b=d$ and $Bc=\eps c$, for $c \in C$.
Under this identification, the triangulated category of mixed
complexes $\DMix$ is equivalent to the derived category of left dg
$\Lambda$-modules.

Let $\cB$ be a small dg-category, i.e., a small category enriched
over graded vector spaces. When $\cB$ is concentrated in degree
$0$ (i.e., when $\cB$ is a $k$-linear category), McCarthy defined
a cyclic module and hence a mixed complex $C_{us}(\cB)$ associated
to $\cB$; see \cite{McC}. Keller observed in \cite[1.3]{Keller99}
that that cyclic module is easily modified to make sense for
general dg-categories. (Since we are working over a field,
Keller's flatness hypothesis is satisfied.)

\begin{Exact-cat}
If $\cA$ is a $k$-linear exact category in the sense of Quillen,
Keller defines the mixed complex $C(\cA)$ in \cite[1.4]{Keller99}
to be the cone of $C_{us}(\Ac^b\cA)\to C_{us}(\bCh^b\cA)$, where
$\bCh^b\cA$ is the dg-category of bounded chain complexes in $\cA$
and $\Ac^b\cA$ is the sub dg-category of acyclic complexes. He
also proves in \cite[1.5]{Keller99} that, up to quasi-isomorphism,
$C(\cA)$ only depends upon the idempotent completion $\cA^+$ of
$\cA$.
\end{Exact-cat}

\begin{ex}\label{ex:CA}
Let $A$ be a $k$-algebra; viewing it as a (dg) category with one
object, $C_{us}(A)$ is the usual mixed complex of $A$ (see
\cite{LodayHC92} or \cite{WeibelHA94}). Now let $\bfP(A)$ denote
the exact category of finitely generated projective $A$-modules.
By McCarthy's theorem \cite[2.4.3]{McC}, the natural map
$C_{us}(A)\to C_{us}(\bfP(A))$ is a quasi-isomorphism of mixed
complexes. Keller proves in \cite[2.4]{Keller99} that
$C_{us}(\bfP(A))\to C(\bfP(A))$ and hence $C_{us}(A)\to
C(\bfP(A))$ is a quasi-isomorphism of mixed complexes. In
particular, it induces quasi-isomorphisms of $HC$, $HP$ and $HN$
complexes.
\end{ex}
\smallskip

\begin{Exactdg}
Let $\cB$ be a small dg-category, and let $DG(\cB)$ denote the
category of left dg $\cB$-modules. There is a Yoneda embedding
$Y:Z^0\cB \to DG(\cB)$, $Y(B)(A)={\cB}(A,B)$, where $Z^0\cB$ is
the subcategory of $\cB$ whose morphisms from $A$ to $B$ are
$Z^0{\cB}(A,B)$. Following Keller \cite[2.1]{Keller99}, we say
that a dg-category is {\it exact} if $Z^0{\cB}$ (the full
subcategory of representable modules $Y(B)$) is closed under
extensions and the shift functor in $DG(\cB)$. The triangulated
category $\cT(\cB)$ of an exact dg-category $\cB$ is defined to be
Keller's stable category $Z^0\cB/B^0\cB$.
\end{Exactdg}

\smallskip
\begin{locpairs}
A {\it localization pair} $\cB=(\cB_1, \cB_0)$ is an exact
dg-category $\cB_1$ endowed with a full dg-subcategory
$\cB_0\subset\cB_1$ such that $Z^0\cB_0$ is an exact subcategory of
$Z^0\cB_1$ closed under shifts and extensions. For a localization
pair $\cB$, the induced functor on associated triangulated
categories $\cT(\cB_0) \subset \cT(\cB_1)$ is fully faithful, and
the {\it associated triangulated category} $\cT(\cB)$ of $\cB$ is
defined to be the Verdier quotient $\cT(\cB_1)/\cT(\cB_0)$.
\end{locpairs}

\begin{sqlocpairs}
\label{M:subLocPair} Let $\cB=(\cB_1, \cB_0)$ be a localization
pair, and let $\cS \subset \cT(\cB)$ be a full triangulated
subcategory. Let $\cC \subset \cB_1$ be the full dg subcategory
whose objects are isomorphic in $\cT(\cB)$ to objects of $\cS$.
Then $\cB_0 \subset \cC$ and $\cC \subset \cB_1$ are localization
pairs, and the sequence $(\cC,\cB_0) \to \cB \to (\cB_1,\cC)$ has
an associated sequence of triangulated categories which is
naturally equivalent to the exact sequence of triangulated
categories $\cS \to \cT(\cB) \to \cT(\cB)/\cS$.
\end{sqlocpairs}

A dg category $\cB$ over a ring $R$ is said to be {\it flat} if
each $H={\cB}(A,B)$ is flat in the sense that $H\otimes_R-$
preserves quasi-isomophisms of graded $R$-modules. A localization
pair $\cB$ is flat if $\cB_1$ (and hence $\cB_2$) is flat. When
the ground ring is a field, as it is in this article, every
localization pair is flat.

In \cite[2.4]{Keller99}, Keller associates to a flat localization
pair $\cB$ a mixed complex $C(\cB)$, the cone of $C(\cB_0)\to
C(\cB_1)$, and proves the following in \cite[Theorem 2.4]{Keller99}:

\smallskip
\begin{thm} \label{KellerLoc}
Let $\cA\to \cB \to \cC$ be a sequence of localization pairs such
that the associated sequence of triangulated categories is exact
up to factors. Then the induced sequence $C(\cA) \to C(\cB) \to
C(\cC)$ of mixed complexes extends to a canonical distinguished
triangle in $\DMix$
$$C(\cA) \to C(\cB) \to C(\cC) \to C(\cA)[1].$$
\end{thm}
\smallskip

\begin{ex} \label{ex:HCX}
The category $\bCh_{\mathrm{parf}}(X)$ of perfect complexes on $X$
is an exact dg-category if we ignore cardinality issues. We need a
more precise choice for the category of perfect complexes. Let $F$
be a field of characteristic zero containing $k$. For $X \in
\SchF$, we choose $\bCh_{\mathrm{parf}}(X)$ to be the category of
perfect bounded above complexes (under cohomological indexing) of
flat $\cO_X$-modules whose stalks have cardinality at most the
cardinality of $F$. (Since $F$ is infinite, all algebras
essentially of finite type over $F$ have cardinality at most the
cardinality of $F$). This is an exact dg-category over $k$. Let
$f:X \to Y$ be a map of schemes essentially of finite type over
$F$. Then $Lf^*$ is $f^*$ on $\bCh_{\mathrm{parf}}(X)$, so
$\bCh_{\mathrm{parf}}$ is functorial up to (unique) natural
isomorphism of functors on $\SchF$. If we want to get a real
presheaf of dg categories on $\SchF$, we can replace
$\bCh_{\mathrm{parf}}$ by the Grothendieck construction
$\int_{\SchF}\bCh_{\mathrm{parf}}$. (This trick is called the {\it
Kleisli rectification}.)

Let $\Ac(X)\subset \bCh_{\mathrm{parf}}(X)$ be the full
dg-subcategory of acyclic complexes. Then
${\bCh}_{\mathrm{parf}}(X)= (\bCh_{\mathrm{parf}}(X), \Ac(X))$ is
a localization pair over $k$ whose associated triangulated
category is naturally equivalent to $\D(X)$ (\cite[3.5.3]{TT},
except for the cardinality part). We define $\bC(X)$ to be the
mixed complex (over $k$) associated to
${\bCh}_{\mathrm{parf}}(X)$.

We define $\bHC(X)$, $\bHP(X)$, $\bHN(X)$ to be the cyclic,
periodic cyclic, negative cyclic homology complexes associated
with the mixed complex $\bC(X)$. In particular, $\bHC$, $\bHP$ and
$\bHN$ are presheaves of complexes on $\SchF$. Keller proves in
\cite[5.2]{KellerDoc} that these definitions agree with the
definitions in \cite{WeibelHC}, with $HC_n(X)=H^{-n}\bHC(X)$, etc.
In addition, the Hochschild homology of $X$ is the homology of the
complex underlying $\bC(X)$.
\end{ex}
\smallskip

\begin{ex}\label{ex:XonZ}
If $Z\subset X$ is closed, let ${\bCh}_{\mathrm{parf}}(X \text{ on
} Z)$ be the localization pair formed by the category of perfect
complexes on $X$ which are acyclic on $X-Z$, and its full
subcategory of acyclic complexes. We define $\bC(X \text{ on } Z)$
to be the mixed complex associated to this localization pair.

If $U\subset X$ is the open complement of $Z$, then Thomason and
Trobaugh proved in \cite[\S5]{TT} that the sequence
${\bCh}_{\mathrm{parf}}(X \text{ on } Z) \to
{\bCh}_{\mathrm{parf}}(X) \to{\bCh}_{\mathrm{parf}}(U)$ is such
that the associated sequence of triangulated categories is exact
up to factors. As pointed out in \cite[5.5]{KellerDoc}, Keller's
Theorem \ref{KellerLoc} implies that $\bC(X \text{ on } Z) \to
\bC(X) \to \bC(U)$ fits into a distinguished triangle in $\DMix$.
\end{ex}

As a consequence of \ref{ex:HCX} and \ref{ex:XonZ}, and a standard
argument involving \'etale covers, we recover the following
theorem, which was originally proven by Geller and Weibel in
\cite[4.2.1 and 4.8]{WG}. (The term ``\'etale descent'' used in
\cite{WG} implies Nisnevich descent; for presheaves of
$\Q$-modules, they are equivalent notions.)

\begin{thm} \label{thm:nis}
Hochschild, cyclic, periodic and negative cyclic homology satisfy
Nisnevich descent.
\end{thm}

We are now ready to prove the cyclic homology analogue of
Thomason's theorem for regular embeddings.
\goodbreak

\begin{thm} \label{thm:HCregdesc}
Let $Y \subset X$ be a regular embedding of $F$-schemes of pure
codimension $d$, let $X' \to X$ the blow-up of $X$ along $Y$ and
$Y'$ the exceptional divisor. Then the presheaves of cyclic,
periodic cyclic and negative cyclic homology complexes satisfy the
Mayer-Vietoris property for the square
$$
\begin{CD}
Y' @>>> X' \\
@VVV @VVV \\
Y @>>> X.
\end{CD}
$$
\end{thm}

\begin{proof}
By section \ref{M:subLocPair}, the filtrations in Proposition
\ref{M:prop:filt} induce filtrations on both
$\bCh_{\mathrm{parf}}(X')$ and on $\bCh_{\mathrm{parf}}(Y')$, and
$Lf^*=f^*$ is compatible with these filtrations. Moreover, $f^*$
induces a map on associated graded localization pairs. By Theorem
\ref{KellerLoc} and Proposition \ref{M:prop:filt}, each square in
the map of filtrations induces a homotopy cartesian square of
mixed complexes, hence the outer square is homotopy cartesian,
too.
\end{proof}

\begin{rem}\label{rem:split}
The filtrations in Proposition \ref{M:prop:filt} split (see proof
of \ref{M:prop:filt}), and induce the usual projective space
bundle and blow-up formulas
$$\bHC(Y')= \bHC(\bP^{d-1}_Y) \simeq \bigoplus_{0 \leq l \leq d-1}
\hskip-9pt\bHC(Y), \ \ {\rm and}\ \ \ \bHC(X')\simeq \bHC(X)\oplus
\bigoplus_{1 \leq l \leq d-1}\hskip-3pt \bHC(Y).
$$
Similarly for $\bHP$ and $\bHN$ in place of $\bHC$. For more
details in the $K$-theory case, see \cite{ThomRD}.
\end{rem}

\begin{rem}\label{rem:fundthm}
Combining the Mayer-Vietoris property for the usual covering of
$X\times\mathbb{P}^1$ with the decomposition of \ref{rem:split}
yields the Fundamental Theorem for negative cyclic homology, which
states that there is a short exact sequence,
$$
0\to \bHN(X\times\A^1)\cup_{\bHN(X)} \bHN(X\times\A^1)  \to
\bHN(X\times(\A^1-\{0\})) \to \bHN(X)[1] \to 0.
$$
This sequence is split up to homotopy; the splitting
$\bHN(X)[1]\to\bHN(X\times(\A^1-\{0\}))$ is multiplication by the
class of $dt/t \in HN_1(k[t,1/t])$. The same argument shows that
there are similar Fundamental Theorems for cyclic and periodic
cylic homology.
\end{rem}

\goodbreak%____________________________________________________
\section{Descent for the $cdh$-topology}

We recall the definition of a $cd$-structure, given in
\cite{VHomotopycd} and \cite{VUnst}.

\begin{defn}\label{def:cd}
Let $\cC$ be a small category. A {\em $cd$-structure} on
$\cC$ is a class $\mathcal{P}$ of commutative squares in
$\cC$ that is closed under isomorphism.
\end{defn}

A $cd$-structure defines a topology on $\cC$. We use the following
$cd$-structures on $\SchF$ and on the subcategory $\SmF$ of
essentially smooth schemes (that is, localizations of smooth
schemes) over $F$.

\begin{ex}\label{ex:cd}
\begin{enumerate}
\item The {\it combined $cd$-structure} on the category $\SchF$.
This consists of all elementary Nisnevich and abstract blow-up
squares. It is complete (\cite[Lemma 2.2]{VUnst}), bounded
(\cite[Proposition 2.12]{VUnst}) and regular (\cite[Lemma
2.13]{VUnst}). By definition,  the $cdh$-topology is the topology
generated by the combined $cd$-structure (see \cite[Proposition
2.16]{VUnst}). \item The  combined $cd$-structure on $\SmF$ is the
sum of the ``upper'' and ``smooth blow-up'' $cd$-structures on
$\SmF.$ It consists of all elementary Nisnevich squares and those
abstract blow-up squares of smooth schemes isomorphic to a blow-up
of a smooth scheme along a smooth center (this $cd$-structure is
discussed in \cite[Section 4]{VUnst}). This $cd$-structure is
complete, bounded and regular (because resolution of singularities
holds over $F$, see the discussion following \cite[Lemma
4.5]{VUnst}). By definition,  the $scdh$-topology is the topology
generated by this $cd$-structure. It coincides with the
restriction of the $cdh$-topology to $\SmF$ (see \cite[Section
5]{SVBK} for more on $cdh$- and $scdh$-topology and their
relationship).
\end{enumerate}
\end{ex}

We shall be concerned with two notions of weak equivalence for a
morphism $f:\EE\to \EE'$ between presheaves of spectra (or
simplicial presheaves) on a category $\cC$. We say that $f$ is a
{\it global weak equivalence} if $\EE(U)\to \EE'(U)$ is a weak
equivalence for each object $U$. If $\cC$ is a site, we say that
$f$ is a {\it local weak equivalence} if it induces an isomorphism
on sheaves of stable homotopy groups (or ordinary homotopy groups,
in the case of simplicial presheaves).

We are primarily interested in the following model structures on
the categories of presheaves of spectra (or simplicial presheaves)
on a category $\cC$; the terminology is taken from \cite{Blander}.
First, there is the {\em global projective model structure} for
global weak equivalences. A morphism $f:\EE \to \EE'$ is a
fibration in this global projective model structure provided
$f(U): \EE(U) \to \EE'(U)$ is a fibration of spectra for each
object $U$ of $\cC$ (we say that weak equivalences and fibrations
are defined objectwise); cofibrations are defined by the left
lifting property. If $\EE\to\EE'$ is a cofibration then each
$\EE(U)\to\EE'(U)$ is a cofibration of spectra, but the converse
does not hold.

Second, for a site $\cC$ there is the {\em local
injective model structure} for local weak equivalences. A morphism
$\EE\to\EE'$ is a cofibration in this model structure if each
$\EE(U)\to\EE'(U)$ is a cofibration; fibrations are defined by the right
lifting property. These model structures were studied by Jardine
in \cite{JardineSPS} and \cite{JardineGen}.

Third, there is the {\em local projective} (or {\em
Brown-Gersten}) model structure for local weak equivalences. A
morphism $\EE \to \EE'$ is a cofibration in this model structure
if it is a global projective cofibration;
fibrations are defined by the right lifting property.

We warn the reader that our local projective model structure for
presheaves is slightly different from (but Quillen-equivalent to)
the corresponding model structure for sheaves discussed in
\cite{VHomotopycd}.

Note that since a cofibration in the local projective model
structure is an objectwise cofibration, it is also a cofibration
in the local injective model structure. In particular, trivial
cofibrations in the local projective structure are also trivial
cofibrations in the local injective model structure. It follows
from the lifting property that a morphism of presheaves
which is a fibration in the local injective model structure is also
a fibration in the local projective model structure.

Any local weak equivalence $\EE\to \EE'$ between local projective
fibrant presheaves is a global weak equivalence. This useful
remark follows from the fact that the identity functor on the
category of presheaves of spectra (respectively, simplicial sets)
is a right Quillen functor from the local projective to the global
projective model structure and hence preserves weak equivalences
between fibrant objects, see \cite[Proposition 8.5.7]{Hirschhorn}.

Recall that a {\it fibrant replacement} of $\EE$ in a model
category is a trivial cofibration $\EE \to \EE'$ with $\EE'$
fibrant. Even though we don't need it, we remark that for all the
model structures we consider, a fibrant replacement can be chosen
functorially by the ``small object argument'' (see
\cite{JardineSPS} for the local injective model structure, and
\cite{Blander} for the projective model structures). We will fix a
fibrant replacement functor $\EE \to \bbH_{\cC}(-,\EE)$ for the
local injective model structure, and we will drop the site from
the notation when the topology is clear from the context.
Following Thomason \cite[p.~532]{AKTEC}, we write $\bbH^n(X,\EE)$
for $\pi_{-n}\bbH(X,\EE)$.

\begin{defn}\label{defn:quasif}
A presheaf of spectra (or simplicial sets) $\EE$ on a site $\cC$
is called {\em quasifibrant} if the local injective fibrant
replacement $\EE\to \bbH(-,\EE)$ is a global weak equivalence,
i.e., the map $\EE(U)\to \bbH(U,\EE)$ is a weak equivalence for
all $U$ in $\cC$.
\end{defn}

An important result of \cite{VHomotopycd} is that under certain
conditions, presheaves satisfying the Mayer-Vietoris property are
precisely quasifibrant presheaves. Note that in that paper,
presheaves satisfying the MV-property are called ``flasque'' .

\begin{thm} \label{thm:BG}
Let $\cC$ be a category with a complete bounded regular
$cd$-structure $\cP$. Then a presheaf of simplicial sets
(or spectra) $\EE$ on $\cC$ is quasifibrant (with respect to
the topology induced by $\cP$) if and only if $\EE$
satisfies the MV-property for $\cP$.
\end{thm}

\begin{proof}
We first prove this for presheaves of simplicial sets. Both
properties, satifying the MV-property and being quasifibrant, are
invariant under global weak equivalence. A local injective fibrant
presheaf is globally equivalent to a fibrant replacement (as
sheaf) of its sheafification. Hence \cite[Lemma 4.3]{VHomotopycd}
shows that a quasifibrant presheaf satisfies the MV-property.
Conversely, \cite[Lemma 3.5]{VHomotopycd} asserts that a local
weak equivalence between presheaves satisfying the MV-property is
a global weak equivalence. As fibrant presheaves satisfy the
MV-property, this implies that any local injective fibrant
replacement $\EE\to \tilde{\EE}$ is a global weak equivalence if the
presheaf satisfies the MV-property; that is, presheaves that
satisfy the MV-property are quasifibrant.

The assertion for presheaves of spectra follows from this, because
a fibrant spectrum is an $\Omega$-spectrum. Indeed, since the
properties ``quasifibrant'' and ``satisfying the MV-property'' are
once again invariant under global weak equivalence, we can as well
assume that all our presheaves are global projective fibrant, in
particular, they are presheaves of $\Omega$-spectra. Now a map of
$\Omega$-spectra is a (stable) weak equivalence if and only if it
is levelwise a weak equivalence of simplicial sets, and a square
of $\Omega$-spectra is homotopy cartesian if and only if it is
levelwise a homotopy cartesian square of simplicial sets. This
reduces the proof to the case of presheaves of simplical sets, as
claimed.
\end{proof}

\begin{terminology}\label{term:descent}
If a presheaf $\EE$ satisfies the equivalent conditions in Theorem
\ref{thm:BG} for a topology $t$ generated by a complete regular
bounded $cd$-structure $\cP$ we say that $\EE$ satisfies
{\em $t$-descent}, or {\em descent for the $t$-topology}.
\end{terminology}

For later use, we note that the analogues of Theorem \ref{thm:BG}
hold for complexes of (pre)sheaves of abelian groups.

\begin{defn}\label{def:compquasif}
Let $\cC$ be a category with a $cd$-structure $\cP$. Let $A^\bu$
be a presheaf of cochain complexes on $\cC$. We say that $A^\bu$
is {\em quasifibrant} for the topology generated by $\cP$ provided
the natural map $A^\bu(U) \to \dR\Gamma(U,A^\bu)$ is a
quasi-isomorphism for each object $U$ of $\cC$. (This property is
usually called ``pseudoflasque'' because it is satisfied by any
cochain complex of flasque sheaves.)

We say that $A^\bu$ satisfies the MV-property for $\cP$, if,
for any square $Q\in \cP$, the square of complexes
$A^\bu(Q)$ is homotopy cocartesian.
\end{defn}

The notation is explained by the following ``great
enlightenment,'' due to Thomason \cite[5.32]{AKTEC}. If $\EE$ is
the presheaf of Eilenberg-Mac~Lane spectra associated to $A^\bu$,
then $\bbH(-,\EE)$ is the presheaf of Eilenberg-Mac~Lane spectra
associated to $R\Gamma(-,A^\bu)$, and we have $\bbH^n(X,\EE)\cong
H^n(X,A^\bu)$.

With these definitions, the exact analogue of Theorem \ref{thm:BG}
holds for complexes of presheaves.

\begin{thm}\label{thm:compBG}
Suppose that $\cC$ is a category with a complete bounded
regular $cd$-structure $\cP$. Then a complex of presheaves
$A^\bu$ is quasifibrant if and only if it satisfies the
MV-property for $\cP$.
\end{thm}

\begin{proof}
Reduce to the result for presheaves of spectra by associating
Eilenberg-Mac~Lane spectra to all complexes.
\end{proof}

\begin{terminology}\label{term:compdescent}
Once again, we say that a presheaf $A^\bu$ of complexes satisfies
$t$-descent for a topology $t$ generated by a $cd$-structure
$\cP$ if it satisfies the equivalent properties of Theorem
\ref{thm:compBG}.
\end{terminology}

\begin{cor}\label{cor:scdh}
Let $A$ be a presheaf of spectra (respectively, complexes) on $\SmF$.
Then $A$ satisfies $scdh$-descent (\ref{ex:cd}.2)
if and only if $A$ satisfies Nisnevich descent and
$A$ satisfies the Mayer-Vietoris property for smooth blow-up squares.
\end{cor}

\begin{ex}\label{ex:singcoh}
It follows from \cite[exp.~Vbis, Corollaire 4.1.6]{SGA4II} that
singular cohomology satisfies $cdh$-descent on the category
$\mathrm{Sch}/\C$. By this we mean that the presheaf of complexes
$X\mapsto S^*(X^{an})$ assigning to a complex variety $X$ the
singular cochain complex of its associated analytic space
satisfies $cdh$-descent.
\end{ex}

For any presheaf $A$ on $\SchF$ write $rA$ for the presheaf on
the subcategory $\SmF$ of smooth schemes obtained by restriction.
The following lemma is immediate from
the observation that $r(a_{cdh}\pi_* A) = a_{scdh}\pi_*(rA)$.

\begin{lem}\label{lem:(s)cdh}
Let $f:A\to B$ be a morphism of presheaves of spectra on $\SchF$.
If $f$ is a local weak equivalence in the $cdh$-topology then $rf:
rA \to rB$ is a local weak equivalence in the $scdh$-topology.
\end{lem}

We can now prove the main technical result of this section; it
will be used in \ref{thm:Kinfdesc} to show that infinitesimal
$K$-theory satisfies $cdh$-descent. Recall that the combined
$cd$-structures on schemes and smooth schemes are complete,
bounded and regular, so that Theorem \ref{thm:BG} applies.

\begin{thm}\label{thm:cdhdesc}
Let $\EE$ be a presheaf of spectra on $\mathrm{Sch}/F$ such that
$\EE$ satisfies excision, is invariant under infinitesimal
extension, satisfies Nisnevich descent and satisfies the
Mayer-Vietoris property for every blow-up along a regular
sequence. Then $\EE$ satisfies $cdh$-descent.
\end{thm}

\begin{proof} We will prove that $\EE$ is $cdh$-quasifibrant.
As $\EE$ satisfies Nisnevich descent and the MV-property for
blow-ups along a regular sequence (in particular, for a blow-up of
a smooth scheme along a smooth subscheme), $r\EE$ satisfies the
MV-property for the combined $cd$-structure on $\SmF$. Let $\EE
\to \bbH_{cdh}(-,\EE)$ be a local injective fibrant replacement of
$\EE$. By Theorem \ref{thm:BG}, $\bbH_{cdh}(-,\EE)$ satisfies the
MV-property for the combined $cd$-structure on $\SchF$. A
fortiori, $r\bbH_{cdh}(-,\EE)$ satisfies the MV-property for the
combined $cd$-structure on $\SmF$. By Lemma \ref{lem:(s)cdh}, the
restriction $r\EE \to r\bbH_{cdh}(-,\EE)$ is a local weak
equivalence in the $scdh$-topology. As source and target satisfy
the MV-property, it is a global weak equivalence on $\SmF$. In
other words, for any smooth scheme $X$, the map $\EE(X) \to
\bbH_{cdh}(X,\EE)$ is a weak equivalence.

Now we proceed as in \cite[Sections 5--6]{HKH}, replacing $\cKH$
by the presheaf $\EE$ everywhere. Specifically, we make the
following conclusions. First of all, because $\EE$ satisfies
excision, Nisnevich descent and is invariant under infinitesimal
extensions, $\EE$ satisfies the MV-property for all finite
abstract blow-ups, as well as for closed covers. If $X$ is a
hypersurface inside some smooth $F$-scheme $U$, we can factor its
resolution of singularities locally into a sequence of blow-ups
along regular sequences and finite abstract blow-ups; using
induction on the dimension of $X$ and the length of the
resolution, we conclude that $\EE(X) \cong \bbH_{cdh}(X,\EE)$ (see
\cite[Theorem 6.1]{HKH} for details of the proof in the case where
$\EE = \cKH$). Next, if $X$ is a local complete intersection, we
use induction on the embedding codimension and Mayer-Vietoris for
closed covers to prove that once again $\EE(X) \cong
\bbH_{cdh}(X,\EE)$ in this case (see \cite[Corollary 6.2]{HKH} for
details). Finally, the general case follows from this because
every integral $F$-scheme is locally a component of a complete
intersection (see \cite[Theorem 6.4]{HKH} for details).
\end{proof}

As a typical application of this result, we prove that periodic
cyclic homology satisfies $cdh$-descent when $\Q\subseteq F$. This
can also be deduced from Feigin-Tsygan's theorem \cite[Theorem
5]{FT} (see also \cite[3.4]{WeibelHo}) \cite[6.8]{CortDR}),  which
identifies HP with crystalline cohomology and from known
properties of the latter established in \cite{Hart}. Note that
$\bHP$ here means the presheaf of complexes computing periodic
cyclic homology over $\Q$.

\begin{cor}\label{cor:HPdesc}
The presheaf of complexes $\bHP$ on $\SchF$ satisfies
$cdh$-descent. Hence its associated presheaf of Eilenberg-Mac~Lane
spectra also satisfies $cdh$-descent.
\end{cor}

\begin{proof}
We have to check the hypotheses of Theorem \ref{thm:cdhdesc} are
satisfied by $\bHP$. The fact that $\bHP$ satisfies excision is in
\cite[5.3]{CQ}; invariance under infinitesimal extensions is
proved in \cite[Theorem II.5.1]{Goodw1}; Nisnevich descent is
Theorem \ref{thm:nis}; and Mayer-Vietoris for blow-ups along a
regular sequence is Theorem \ref{thm:HCregdesc}.
\end{proof}

\begin{ex}\label{ex:KHdesc}
Theorem \ref{thm:cdhdesc} applies in particular to prove that
homotopy $K$-theory $KH$ satisfies $cdh$-descent, as explained in
\cite{HKH}. As a consequence we have the following computation
(see \cite[Theorem 7.1]{HKH}).

Suppose that $X$ is a scheme, essentially of finite type over a
field $F$ of characteristic $0$ and such that $\mathrm{dim}(X) =
d$. Then $KH_n(X) = 0$ for $n < -d$ and $KH_{-d}(X) =
H^d_{cdh}(X,\Z)$.
\end{ex}

\section{Descent for infinitesimal $K$-theory}

In this section, we combine the previous sections to prove (in
Theorem \ref{thm:Kinfdesc} below) that Corti\~nas' infinitesimal
$K$-theory satisfies $cdh$-descent on $\mathrm{Sch}/F$, for any
field $F$ of characteristic $0$. All variants of cyclic homology
are taken over $\Q$.

Recall from \ref{ex:HCX} that $\bHN(X)$ is the presheaf of
complexes defining negative cyclic homology; we obtain a presheaf
of spectra from this by taking the associated Eilenberg-Mac~Lane
spectrum. There is a Chern character $\cK(X) \to \mathbf{HN}(X)$
(for a definition, start for example with \cite[Section
4]{WeibelNil}, use the Fundamental Theorem \ref{rem:fundthm} to
extend to non-connective $K$-theory and globalize using Zariski
descent). Here $\cK(X)$ denotes the non-connective $K$-theory
spectrum of perfect complexes on $X$ as in \cite[6.4]{TT}.

\begin{defn}\label{def:infK}
Let $X$ be a $\Q$-scheme. We define the infinitesimal $K$-theory
of $X$, $\cK^{inf}(X)$, to be the homotopy fiber of the Chern
character $\cK(X) \to \bHN(X)$.
\end{defn}

The following theorem was proven by Corti\~nas in \cite{CortKABI}.
It verified the ``KABI-conjecture'' of
Geller-Reid-Weibel (\cite[0.1]{GRW}).

\begin{thm}[Corti\~nas]\label{thm:KABI}
$\cK^{inf}$ satisfies excision on the category of $\Q$-algebras.
\end{thm}

\begin{thm}\label{thm:KinfNisdesc}
$\cK^{inf}$ satisfies Nisnevich descent.
\end{thm}

\begin{proof}
Both $\cK$ and $\mathbf{HN}$ do, by \cite[10.8]{TT} and Theorem
\ref{thm:nis}.
\end{proof}

From Theorem \ref{thm:HCregdesc} and Remark \ref{M:rem2}, both
$\bHN$ and $K$-theory have the Mayer-Vietoris property for any
square associated to a blow-up along a regular embedding. This
proves the following result.

\begin{thm}\label{thm:Kinfregdesc}
$\cK^{inf}$ satisfies the Mayer-Vietoris property for every
blow-up along a regular embedding.
\end{thm}

Finally, we have the following result, due to Goodwillie, see
\cite{Goodw2}.

\begin{thm}[Goodwillie]\label{thm:infext}
Let $A$ be a $\Q$-algebra and $I \subset A$ a nilpotent ideal.
Then $\cK^{inf}(A,I)$ is contractible. That is, $\cK^{inf}$ is
invariant under infinitesimal extension.
\end{thm}

\begin{proof}
Goodwillie proves (in \cite[Theorem II.3.4]{Goodw2}) that the
Chern character induces an equivalence $\cK(A,I) \to
\mathbf{HN}(A,I)$. This immediately implies the assertion.
\end{proof}

\begin{thm}\label{thm:Kinfdesc}
The presheaf of spectra $\cK^{inf}$ satisfies $cdh$-descent.
\end{thm}

\begin{proof}
This follows from Theorem \ref{thm:cdhdesc}, once we observe that the
presheaf $\cK^{inf}$ satisfies the conditions given in the theorem.
These conditions hold by
Theorem \ref{thm:KABI}, Theorem \ref{thm:Kinfregdesc}, Theorem
\ref{thm:KinfNisdesc} and Theorem \ref{thm:infext}.
\end{proof}

\section{The obstruction to homotopy invariance}

We will say that a sequence of presheaves of spectra
$\EE_1\to\EE_2\to\EE_3$ is an (objectwise) {\em homotopy fibration
sequence} provided that for each scheme $X$, the sequence of
spectra $\EE_1 (X) \to \EE_2(X) \to \EE_3 (X)$ is weakly
equivalent to a fibration sequence (that is, it defines a
distinguished triangle in the homotopy category of spectra). We
have the following useful observation (cf.\ \cite[1.35]{AKTEC},
\cite[p.~73]{JardineSPS}, \cite[p.~194]{JardineLSS}): if
$\EE_1\to\EE_2\to\EE_3$ is a homotopy fibration sequence, then
$\bbH_{cdh}(-,\EE_1)\to\bbH_{cdh}(-,\EE_2)\to\bbH_{cdh}(-,\EE_3)$
is also a homotopy fibration sequence.

For a presheaf of spectra $\EE$ on $\mathrm{Sch}/F$, we will write
$\tilde{C}_j\EE$ for the cofiber of the map
$\EE\to\EE(-\times\A^j)$. Since $\tilde{C}_j\EE$ is a direct
factor of $\EE(-\times\A^j)$, the functor $\EE \mapsto
\tilde{C}_j\EE$ preserves homotopy fibration sequences. If we also
use the $\tilde{C}_j$ notation for presheaves of abelian groups,
then we have $\tilde{C}_j(\pi_r\EE) \cong \pi_r(\tilde{C}_j\EE)$.

\begin{lem}\label{lem:nil1}
Suppose the presheaf of spectra $\EE$ satisfies descent for the
$cdh$-topology (or Zariski, or Nisnevich topology). Then so do the
presheaves $\tilde{C}_j\EE$.
\end{lem}

\begin{proof}
All three topologies are generated by a complete bounded regular
$cd$ structure. By Theorem \ref{thm:BG}, the presheaf $\EE$
satisfies the MV-property, hence so do the presheaves
$\EE(-\times\A^j)$, for all $j$; consequently, the presheaves
$\tilde{C}_j\EE$ also satisfy Mayer-Vietoris, that is, they
satisfy descent.
\end{proof}

In particular, the presheaves $\tilde{C}_j\cK^{inf}$ satisfy
descent for the $cdh$-topology, and the presheaves
$\tilde{C}_j\mathbf{HN}$ satisfy descent for the Zariski topology.

We will say that a presheaf $\EE$ is {\it contractible} if
$\EE(X)\simeq\ast$ for all $X$.

\begin{lem}\label{lem:nil2}
The presheaves $\bbH_{cdh}(-,\tilde{C}_j\cK)$ are contractible
for all $j\ge1$.
\end{lem}

\begin{proof}
As smooth schemes are $K_m$-regular for any $m$, $a_{cdh}\pi_m
\tilde{C}_j\cK = a_{cdh}\tilde{C}_jK_m = 0$ for all $j\ge1$ and all $m$.
The assertion follows from the general local-to-global spectral sequence in
\cite[1.36]{AKTEC}, applied to the $cdh$ site.
\end{proof}

In characteristic zero, we also have the following result, see
\cite{Goodw1}.

\begin{prop}\label{prop:nil3}
For all $j\ge1$, the presheaves $\tilde{C}_j\bHP$ are contractible,
and hence so are the presheaves
$\bbH_{cdh}(-,\tilde{C}_j\bHP)$.
\end{prop}

\begin{proof}
Since $\Q\subseteq F$, $\bHP$ is $\A^1$-homotopy invariant on algebras,
by \cite[III.5.1]{Goodw1} (see also \cite[E.5.1.4]{LodayHC92}).
As $\bHP$ satisfies Zariski descent, this implies the
first assertion. The second assertion is an immediate consequence.
\end{proof}

\begin{cor}\label{cor:nil4}
For all $j\ge1$, there is a (global) weak equivalence
$\tilde{C}_j\mathbf{HC} \cong \Omega \tilde{C}_j\mathbf{HN}$,
and hence an (objectwise) homotopy fibration sequence:
\[
\tilde{C}_j\mathbf{HC} \to \tilde{C}_j\cK^{inf} \to \tilde{C}_j\cK.
\]
\end{cor}

\begin{proof}
Immediate from Proposition \ref{prop:nil3} and the fundamental
homotopy fibration sequences $\bHN(X) \to \bHP(X) \to
\Omega^{-2}\bHC(X)$.
\end{proof}

Applying $\bbH_{cdh}(-,-)$ to this homotopy fibration sequence,
and using Lemma \ref{lem:nil2}, we see that Theorem
\ref{thm:Kinfdesc} implies the next result.

\begin{thm}\label{thm:nilsequence} Let $j \geq 0$.
There is an (objectwise) homotopy fibration sequence
\[ \tilde{C}_j\mathbf{HC} \to \bbH_{cdh}(-,\tilde{C}_j\mathbf{HC}) \to
\tilde{C}_j\cK \]
\end{thm}

\begin{lem}\label{lem:negNHC} Let $j\geq 0$.
The Zariski sheaves $a_{Zar}\pi_n \tilde{C}_j\mathbf{HC}$ (and a
fortiori, the $cdh$-sheaves $a_{cdh}\pi_n
\tilde{C}_j\mathbf{HC}$) vanish for all $n < 0$.
\end{lem}

\begin{proof}
For any ring $A$, $\bHC(A)$ is $-1$-connected. Therefore
$\tilde{C}_jHC_n(A) = 0$ for $n < 0$. This implies the assertion.
\end{proof}

\begin{rem}
The vanishing range in \ref{lem:negNHC} is best possible,
because $a_{Zar}\pi_0\tilde{C}_j\bHC$ is $\tilde{C}_j\cO$.
\end{rem}

\begin{cor}\label{cor:HdX}
For all $m>d=\dim(X)$ and all $j$,
$HC_{-m}(X\times\A^j) = \bbH^m_{cdh}(X\times\A^j,\bHC) =0$.
Moreover, $\bbH^d_{Zar}(X,\tilde{C}_j\bHC)=H^d_{Zar}(X,\tilde{C}_j\cO)$
and $\bbH^d_{cdh}(X,\tilde{C}_j\bHC)=H^d_{cdh}(X,a_{cdh}\tilde{C}_j\cO)$.
\end{cor}

\begin{proof}
Since both the Zariski and $cdh$ cohomological dimension of $X$
are at most $d$, these follow from \ref{lem:negNHC} via the Leray
spectral sequences $H^p(X,a_{t}\pi_{-q}\tilde{C}_j\bHC)\Rightarrow
H^{p+q}(X,\tilde{C}_j\bHC)$.
\end{proof}

From this we can already conclude the following weak form of
Conjecture \ref{conj:reg}.

\begin{cor}\label{cor:Kdimweak}
Let $F$ be a field of characteristic $0$, and $X$ a
$d$-dimensional scheme, essentially of finite type over $F$. Then
$X$ is $K_{n}$-regular and $K_n(X) = 0$ for all $n < -d$.
\end{cor}

\begin{proof}
The first part is an immediate consequence of Theorem
\ref{thm:nilsequence} and Corollary \ref{cor:HdX}. The second part
follows from the first using the spectral sequence
$K_q(X\times\A^p) \Rightarrow KH_{p+q}(X)$.
\end{proof}

\begin{cor}\label{cor:Hd-seq}
If $\dim(X)=d$, there is an exact sequence for every $j\ge1$:
\[
H^d_{Zar}(X,a_{Zar}\tilde{C}_jHC_0) \to
H^d_{cdh}(X,a_{cdh}\tilde{C}_jHC_0)\to \tilde{C}_jK_{-d}(X) \to 0.
\]
\end{cor}

\begin{proof}
Combine \ref{thm:nilsequence} and \ref{cor:HdX}.
\end{proof}

\section{$cdh$-cohomology of coherent sheaves and the K-dimension
conjecture}

Most of this section will be taken up by the proof of the next
result.

\begin{thm}\label{thm:topchange} Let $X$ a $F$-scheme, essentially
of finite type, and of dimension $d$. Then the natural
homomorphism, induced by the change of topology,
\[ H^d_{Zar}(X,\cO_X) \longrightarrow H^d_{cdh}(X,a_{cdh}\cO_X) \]
is surjective.
\end{thm}

Before proving this theorem, we show how Theorem
\ref{thm:topchange} and Corollary \ref{cor:Hd-seq} imply the
K-dimension Conjecture \ref{conj:reg}.

\begin{thm}\label{thm:Kdim}
Let $F$ be a field of characteristic $0$ and $X$ be an $F$-scheme,
essentially of finite type and of dimension $d$. Then $X$ is
$K_{-d}$-regular and $K_n (X) = 0$ for $n < -d$. Moreover,
$K_{-d}(X) \cong H^d_{cdh}(X,\Z).$
\end{thm}

\begin{proof} Fix $j\ge1$ and let $V_j$ denote the $F$-vector space
$F[t_1,\dotsc t_j]/F$.
Since $HC_0(A) = A$ for any commutative algebra $A$,
$\tilde{C}_jHC_0(A) = A\otimes_F V_j$.
Hence $a_{Zar}\tilde{C}_jHC_0 \cong \cO_X\otimes_F V_j$ and
$a_{cdh}\tilde{C}_jHC_0 \cong a_{cdh}\cO_X\otimes_F V_j.$
Therefore Theorem \ref{thm:topchange} and Corollary \ref{cor:Hd-seq}
imply that $X$ is $K_{-d}$-regular.

The remaining assertions follow from the calculation of $KH_*(X)$
in Example \ref{ex:KHdesc}, and the spectral sequence
$K_q(X\times\A^p) \Rightarrow KH_{p+q}(X)$.
\end{proof}

The proof of Theorem \ref{thm:topchange} will be in two parts:
First, we will prove a stronger result for smooth $X$. The second
part is the proof for a general $X$. To simplify notation, for any
topology $t$ and presheaf $A$, we write $H^*_t(X,A)$ for
$H^*_t(X,a_t A)$. We write $\cO$ for the presheaf $X\mapsto
\cO_X(X)$, and $\dR\Gamma_t(X,\cO)$ for
a functorial model for the total right derived functor of the
global sections functor $X\mapsto a_t\cO(X)$.

\begin{prop}\label{prop:smoothtc}
Let $X$ be a smooth $F$-scheme. Then $\cO(X) \cong a_{cdh}\cO(X)$
and the natural homomorphism
\[ H^*_{Zar}(X,\cO) \longrightarrow H^*_{cdh}(X,\cO)\]
is an isomorphism.
\end{prop}

\begin{proof} We need to show that the natural map
$\dR\Gamma_{Zar}(X,\cO)\to \dR\Gamma_{cdh}(X,\cO)$ is a
quasi-isomorphism. Since the target satisfies $cdh$-descent
\ref{term:compdescent}, this amounts to showing that the presheaf
of complexes $X\mapsto \dR\Gamma_{Zar}(X,\cO)$ on $\SmF$ satisfies
the conditions of Corollary \ref{cor:scdh}. First of all, it is
classical that $\dR\Gamma_{Zar}(X,\cO) \cong
\dR\Gamma_{Nis}(X,\cO)$, which implies that this presheaf
satisfies Nisnevich descent (it sends elementary Nisnevich squares
to homotopy cocartesian squares). Hence, it suffices to show that
$\dR\Gamma_{Zar}(-,\cO)$ transforms smooth blow-up squares into
homotopy cocartesian squares. To this end, let $X$ be a smooth
scheme, $Y\subset X$ a smooth closed subscheme, $p:X' \to X$ the
blow-up along $Y$ and $j:Y'\subset X'$ the exceptional divisor;
write $q: Y' \to Y$ for the restriction. Then we need to show that
the natural map
\[
\mathrm{Cone}\biggl(\dR\Gamma_{Zar}(X,\cO) \to\dR\Gamma_{Zar}(X',\cO)\biggr)
 \longrightarrow
\mathrm{Cone}\biggl(\dR\Gamma_{Zar}(Y,\cO) \to\dR\Gamma_{Zar}(Y',\cO)\biggr)
\]
is a quasi-isomorphism.
In fact, both of those cones are $0$. Indeed,
$\cO_X \to \dRp_*p^*\cO_X$ is a quasi-isomorphism by
\cite[Lemme 2.3(a)]{ThomRD}, and
the usual computation of cohomology of projective space shows that
$\cO_Y \to \dRq_*\dLq^*\cO_Y = \dRq_*\cO_{Y'}$
is also a quasi-isomorphism.
Hence we have a homotopy cartesian square, or alternatively,
\[ \dR\Gamma_{Zar}(X,\cO) \to
\dR\Gamma_{Zar}(X',\cO) \times\dR\Gamma_{Zar}(Y,\cO) \to
\dR\Gamma_{Zar}(Y',\cO)
\]
is a homotopy fibration sequence. It follows that
$\dR\Gamma_{Zar}(-,\cO)$ satisfies $scdh$ descent, and in
particular that $\cO\cong a_{cdh}\cO$ on $\SmF$.
\end{proof}

Neither the global sections of $\cO$ nor the higher Zariski
cohomology of $\cO$ satisfies $cdh$-descent for non-smooth
schemes; this fails for example when $X$ is a cusp. Nevertheless,
we have the following partial result, which suffices for Theorem
\ref{thm:topchange}.

\begin{lem}\label{lem:zarepi}
Let $X$ be a reduced affine Noetherian scheme and $p: X'\to X$ a
proper morphism such that all the fibers of $p$ have dimension at
most $d-1$. Let $j:Y'\subset X'$ be a closed subscheme. Then the
restriction map $H^{d-1}_{Zar}(X',\cO_{X'}) \to
H^{d-1}_{Zar}(Y',\cO_{Y'})$ is surjective.
\end{lem}

\begin{proof}
The theorem on formal functions implies that $R^d p_*\cF = 0$ for
any quasicoherent sheaf $\cF$ on $X'$. Hence the functor
$R^{d-1}p_*$ is right exact on quasicoherent sheaves. In
particular, $R^{d-1}p_*\cO_{X'}\to R^{d-1}p_*j_*\cO_{Y'}$ is onto.
Because $X$ is affine, this proves the assertion.
\end{proof}

For legibility, we will write $a$ for the natural morphism of sites
from the $cdh$-site to the Zariski site on $\SchF$. If $\cF$ is a Zariski
sheaf, then $a^*\cF$ is the same as $a_{cdh}\cF$.

\begin{lem}\label{lem:coh} For any scheme $X$ of finite type over
$F$, the complex of Zariski sheaves $\dR a_* a^*\cO\vert_{X_{Zar}}$
has coherent cohomology sheaves.
\end{lem}

\begin{proof} If $X$ is smooth then the assertion is an immediate
consequence of Proposition \ref{prop:smoothtc}. We prove the
general case by induction on the dimension of $X$. If $\dim(X)=0$,
then $X=\Spec(A)$ for some Artinian ring $A$ and $\dR a_* a^*\cO
\cong \cO^{red}$, which is a coherent sheaf. Now suppose
$d=\dim(X)> 0$. Let $p:X' \to X$ be a resolution of singularities,
$i:Y \subset X$ the singular set and $j:Y'\subset X'$ the
exceptional divisor. Because $Ra_*$ commutes with $Rf_*$ for every
morphism $f$ in $Sch/F$, we have a distinguished triangle of
complexes of sheaves of $\cO_X$-modules on $X_{Zar}$:
\[
\dR a_* a^*\cO \to \dRp_*\dR a_* a^*\cO \times \dR i_*\dR a_* a^*\cO
\to \dR(pj)_*\dR a_* a^*\cO.
\]
The second and third terms in this triangle have coherent
cohomology sheaves; this follows from induction on the dimension,
the assertion for the smooth $X'$, and the fact that proper
morphisms have coherent direct images. Hence, the first term has
coherent cohomology sheaves, too (see \cite[Expos\'e I,
Cor.~3.4]{SGA6}).
\end{proof}

\begin{proof}[Proof of Theorem \ref{thm:topchange}]
We proceed by induction on the dimension. If $d = 0$, then
$X =\Spec(A)$ for some Artinian ring $A$, and $H^0(X,\cO) = A
\to A^{red} = H^0_{cdh}(X,\cO)$ is surjective. Now assume that we
have shown the assertion for all schemes of dimension less than $d>0$.

We claim that it suffices to prove the assertion when $X$ is
affine, or indeed for local $X$ of dimension $d$. To see this,
suppose that $X$ is any $d$-dimensional scheme, essentially of
finite type over $F$. We have a Leray spectral sequence
\[
H^p_{Zar}(X,R^q a_* a^*\cO) \Longrightarrow H^{p+q}_{cdh}(X,\cO).
\]
Fix $q>0$ and consider the stalk of $R^qa_* a^*\cO$ at a point
$x\in X$ of codimension $c$ (that is, where the local ring
$\cO_{X,x}$ has Krull dimension $c$). By assumption, the stalk is
zero if $q\ge c$. Since the sheaf $R^q a_* a^* \cO$ is coherent by
Lemma \ref{lem:coh}, this implies that  $R^qa_*a^*\cO$ is
supported on a closed subscheme of codimension $>q$, i.e., of
dimension $< d - q$. This implies that
$H^{p}_{Zar}(X,R^qa_*a^*\cO) = 0$ provided $p+q\ge d$ and $q>0$.
Hence the Leray spectral sequence degenerates enough to show that
$H^d_{Zar}(X,a_*a^*\cO) \to H^d_{cdh}(X,\cO)$ is surjective.

Consider the cokernel $\cF$ of the adjunction map $\cO\to a_*a^*\cO$,
which is coherent by \ref{lem:coh}. It vanishes on an open dense subset
of $X$ (namely, on the complement of the singular set of $X^{red}$,
by \ref{prop:smoothtc}), so $\cF$ is supported in dimension $<d$
and hence $H^d_{Zar}(X,\cF)=0$. Since $H^d_{Zar}(X,-)$ is right exact,
$H^d_{Zar}(X,\cO) \to H^d_{zar}(X,a_* a^* \cO)$ must be a surjection.
This establishes our claim.
To summarize, it suffices to assume the result true in dimension $<d$
and prove it for affine schemes of dimension $d$.
To simplify matters, we can also assume that $X$ is reduced.
Indeed, since $H^d_{Zar}(X,-)$ is right exact, the map
$H^d_{Zar}(X,\cO) \to H^d_{Zar}(X^{red},\cO)$ is surjective, and
$H^d_{cdh}(X,\cO) = H^d_{cdh}(X^{red},\cO)$.

Let $X$ be an affine $d$-dimensional scheme, and choose a
resolution of singularities $p:X' \to X$. Let $Y\subset X$ be the
singular subscheme and $Y'\subset X'$ the exceptional divisor.
Since $Y$ and $Y'$ have smaller dimension, $H^d_{cdh}(Y,\cO) =
H^d_{cdh}(Y',\cO) = 0$ for cohomological dimension reasons
\cite{SVBK}. Furthermore, $p$ is proper and has fibers of
dimension $<d$; because $X$ is affine, this implies that
$H^d_{Zar}(X',\cO) = 0$, by the theorem on formal functions. Since
$X'$ is smooth, we conclude that $H^d_{cdh}(X',\cO) = 0$ by
Proposition \ref{prop:smoothtc}. Now the long exact sequence in
$cdh$-cohomology for the abstract blow-up $p$ gives a diagram with
exact top row.
$$
\begin{CD}
H^{d-1}_{cdh}(Y,\cO)\times
H^{d-1}_{cdh}(X',\cO) @>>> H^{d-1}_{cdh}(Y',\cO) @>>> H^d_{cdh}(X,\cO) @>>>0\\
{\hspace{15ex}}@AAA @AA{\text{onto}}A \\
\hspace{15ex}H^{d-1}_{Zar}(X',\cO) @>{\text{onto}}>> H^{d-1}_{Zar}(Y',\cO)
\end{CD}
$$
The right vertical map in this diagram is
surjective by induction. As the fibers of $p$ have dimension less
than $d$, Lemma \ref{lem:zarepi} implies that the bottom horizontal map
is surjective. Therefore $H^{d-1}_{cdh}(X',\cO) \to H^{d-1}_{cdh}(Y',\cO)$
is also surjective and hence $H^d_{cdh}(X,\cO) = 0$. This finishes
the induction step and the proof of Theorem \ref{thm:topchange}.
\end{proof}

\subsection*{Acknowledgments}
This paper grew out of discussions the authors had at the Institut
Henri Poincar\'e during the semester on $K$-theory and
noncommutative geometry in Spring 2004. We thank the organizers,
M. Karoubi and R. Nest, as well as the IHP for their hospitality.

\end{document}